\documentclass[a4paper,12pt]{article}
\usepackage[norelsize]{algorithm2e}
\usepackage{amssymb}
\usepackage{amsbsy}
\usepackage{amsmath}
\usepackage{graphicx}
\usepackage{arcs}

\newtheorem{theorem}{\bf Theorem}
\newtheorem{definition}{\bf Definition}[section]
\newtheorem{proposition}[definition]{\bf Proposition}
\newtheorem{lemma}[definition]{\bf Lemma}
\newtheorem{corollary}[definition]{\bf Corollary}

\newcommand{\relR}{\mathrel{R}}

\newcommand{\Z}{\mathbb{Z}}
\newcommand{\C}{\mathbb{C}}

\newcommand{\II}{\mathcal{I}}
\newcommand{\F}{\mathcal{F}}
\newcommand{\R}{\mathbb{R}}

\newcommand{\fsf}{f_\mathrm{sf}}
\newcommand{\fzig}{f_\mathrm{zz}}
\newcommand{\fizig}{f_\mathrm{izz}}
\newcommand{\ftent}{f_\mathrm{tent}}

\renewcommand{\Re}{\mathrm{Re}}

\renewcommand{\Im}{\mathrm{Im}}

\newcommand{\ee}{\mathbf{e}}
\newcommand{\ff}{\mathbf{f}}
\renewcommand{\gg}{\mathbf{g}}
\newcommand{\vv}{\mathbf{v}}
\newcommand{\ww}{\mathbf{w}}
\newcommand{\xx}{\mathbf{x}}
\newcommand{\yy}{\mathbf{y}}
\newcommand{\zero}{\mathbf{0}}

\newcommand{\tJ}{\widetilde{J}}
\newcommand{\hP}{\widehat{P}}
\newcommand{\tA}{\widetilde{A}}
\newcommand{\tP}{\tilde{P}}
\newcommand{\tn}{\tau^{(N)}}
\newcommand{\tni}{\tau^{(N_i)}}
\newcommand{\tmix}{\tau_\mathrm{mix}^{(N)}}
\newcommand{\tperm}{\tau_\mathrm{perm}}

\begin{document}
\title{Mixing properties of the zigzag map under composition with permutations, \\ II: Maps of non-constant orientation}

\author{Nigel P.~Byott$^1$, Yiwei Zhang$^2$, Congping Lin$^1$ \\
College of Engineering, Mathematics and Physical Sciences,\\
University of Exeter $^1$\\
Facultad de Matem\'atica\\
Pontificia Universidad Cat\'olica de Chile$^2$\\}

\date{\today}
\maketitle

\begin{abstract}
For an integer $m \geq 2$, let $\mathcal{P}_m$ be the partition of
the unit interval $I$ into $m$ equal subintervals, and let $\F_m$ be
the class of piecewise linear maps on $I$ with constant slope $\pm
m$ on each element of $\mathcal{P}_m$. We investigate the effect on
mixing properties when $f \in \F_m$ is composed with the interval
exchange map given by a permutation $\sigma \in S_N$ interchanging
the $N$ subintervals of $\mathcal{P}_N$. This extends the work in a
previous paper [N.P.~Byott, M.~Holland and Y.~Zhang, DCDS,~{\bf 33},
(2013)~3365--3390], where we considered only the
``stretch-and-fold'' map $\fsf(x)=mx \bmod 1$.
\end{abstract}

\maketitle

\section{Introduction and statement of results}  \label{intro}

A natural question about a dynamical system is how fast it mixes
measurable subsets of its domain.  In the setting of discrete time
piecewise smooth one-dimensional expanding maps, the rate of decay
of correlations for functions of bounded variation gives a
quantitative interpretation of the speed of mixing, and is governed
by the isolated spectrum of the transfer operator on the space of
such functions (see \cite{Baladi00, Boyarsky-Gora} and the
references therein).  As there is no general approach to determine
this isolated spectrum explicitly, detailed studies of specific
families of maps, as for example in \cite{CE04} and \cite{DFS00},
are an important step towards a better understanding of quantitative
mixing phenomena.

For simplicity, we focus our discussion on piecewise linear Markov
maps. A nice property of such maps is that the corresponding
transfer operator has a finite matrix representation. It is perhaps
surprising that, even in this special situation, the precise
determination of the mixing rate is already a non-trivial task
\cite{SBJ13,SBJ13b}. On the other hand, using a Ulam-like
construction \cite{BIS95,Ulam64}, every expanding (non-linear)
Markov map can be approximated by a sequence of piecewise linear
maps. In many circumstances, we can use the mixing rate of high
order piecewise linear maps to bound or approximate the mixing rate
of the original non-linear maps \cite{CPR,SBJ13,SBJ13b}. (This is
extended in \cite{Froy95,Froy97} to the settings of
multi-dimensional expanding maps and of Anosov maps.) Historically,
on the basis of explicit calculations for some topologically mixing
piecewise linear Markov maps (namely the skew doubling maps and skew
tent maps), Badii et al.~\cite{BHMP88} conjectured that, for every
topologically mixing and expanding Markov map, the mixing rate can
be universally bounded from above in terms of a (generalized)
Lyapunov exponent. A counterexample was soon obtained \cite{CPR} by
making a non-linear perturbation to these maps. Indeed, as observed
in \cite{DFS00,Just90,SBJ13,SBJ13b}, the Lyapunov exponent and the
topological entropy only provide a bound from below for the mixing
rate via bounded variation observations.

In a continuous time setting, a novel aspect of mixing was
considered in \cite{Ashwin02}. The domain of a one-dimensional
diffusion process was divided into equal subintervals, and the
underlying function was composed with an interval exchange map
corresponding to a permutation of these. This typically results in
faster mixing than the diffusion alone.

In view of these results, it is of interest to investigate the
analogous effect of composition with a permutation for amenable
examples of piecewise linear Markov maps. Such an investigation
began in our previous paper \cite{BHZ}. In particular, we considered
the ``stretch-and-fold'' map $\fsf(x)=mx \bmod 1$ for an integer $m
\geq 2$ on $I=[0,1]$. (When $m=2$, this is the well-known ``doubling
map''). When $I$ is divided into $N \geq m$ equal subintervals, and
$f$ is composed with a permutation $\sigma$ of these, the Lyapunov
exponent and topological entropy are unchanged, but, in contrast to
the results of \cite{Ashwin02}, mixing typically becomes slower.
Indeed, there may be permutations $\sigma$ such that the new map
$\sigma \circ \fsf$ fails to be topologically mixing. Under the
hypothesis $\gcd(m,N)=1$, we showed that $\sigma \circ \fsf$ is
topologically mixing for all permutations $\sigma$, and we
determined explicitly the worst mixing attained as $\sigma$ varies.
In fact, the mixing can be made arbitrarily slow by allowing $N$ to
increase. This not only shows a striking contrast between the effect
of composing with a permutation in the continuous and discrete time
settings, but also adds to the collection of piecewise expanding
interval maps where mixing rates can be explicitly determined.
Indeed, these maps provide an alternative counterexample to the
conjecture of Badii et al., in which no non-linear perturbations are
required: the Markov structure is changed in an essentially
combinatorial manner, and only piecewise linear maps are used.
Moreover, since we consider Markov partitions with equal
subintervals, this has the advantage that the mixing rate can be
calculated exactly rather than approximated numerically.

In the present paper, we continue this investigation by regarding
$\fsf$ as the simplest member of a family $\F_m$ of interval maps
with piecewise constant integer slope $\pm m$. More precisely, for
fixed $m \geq 2$, and for any sequence $\epsilon_1, \ldots,
\epsilon_m \in \{-1,1\}^m$, we take $f_{\epsilon_1, \ldots ,
\epsilon_m}$ to be the interval map with slope $m\epsilon_j$ on the
interval $((j-1)/m, j/m)$. Explicitly,
$$ f_{\epsilon_1 , \ldots , \epsilon_m}(x) =
     \begin{cases}  mx - j+1 & \mbox{if } j-1 \leq mx < j \mbox{ with } \epsilon_j=1; \cr
                    j - mx & \mbox{if } j-1 \leq mx < j \mbox{ with } \epsilon_j =-1. \cr
             \end{cases}  $$
We consider the family of $2^m$ maps $\F_m=\{ f_{\epsilon_1 , \ldots
,
  \epsilon_m} \, : \, \epsilon_j = \pm 1\}$. Taking $\epsilon_j=1$
for all $j$ gives the stretch-and-fold map $\fsf$ with piecewise
constant slope $m$. At the other extreme, taking
$\epsilon_j=(-1)^{j-1}$ gives a map with $m-1$ changes of slope. We
shall refer to this as the $m$-fold zigzag map, writing $\fzig$ for
$f_{1,-1,1,\ldots}$.  (In the case $m=2$, the zigzag map $f_{1,-1}$
is just the familiar ``tent map'' $\ftent$.) We could equally well
take the inverted zigzag map with $\epsilon_j=(-1)^j$, and we write
$\fizig$ for $f_{-1,+1,-1,\ldots}$.

Let $\mathcal{P}_N$ be the partition of $I$ into $N \geq m$ equal
subintervals. We write $S_N$ for the symmetric group, consisting of
all $N!$ permutations $\sigma$ of $\{1, 2,\ldots,N\}$. Abusing
notation, we also write $\sigma$ for the interval exchange map given
by the corresponding permutation of the subintervals of
$\mathcal{P}_N$:
$$  \sigma(x) = x + (\sigma(j)-j)/N  \mbox{ if } (j-1)/N \leq x <
  j/N.  $$
For each $f \in \F_m$ and each $\sigma \in S_N$, the map $g = \sigma
\circ f$ has piecewise constant slope $\pm m$. It follows that $g$
has topological entropy $\log m$, and that the Lyapunov exponent of
$g$ is also $\log m$; these quantities are therefore unaffected by
composition with $\sigma$ and hence independent of $N$. The purpose
of this paper is to investigate how the mixing rate is affected by
the combinatorial operation of composition with $\sigma$. In
particular, we shall determine to what extent the mixing properties
of $f$ are neutralised by composition with an appropriate choice of
$\sigma$.

Before describing in detail the results of this paper, we explain
the notion of mixing rate we shall use. We quantify the speed of
mixing in terms of the decay of correlations for observables of
bounded variation. The context is briefly summarised below, and we
refer the readers to \cite{Boyarsky-Gora} for more detailed
information. Let $g \colon I \longrightarrow \R$ be a piecewise
expanding map and suppose for the moment that $g$ is topologically
mixing. Then $g$ admits a unique absolutely continuous invariant
probability measure $\mu$. Let $BV$ denote the Banach space of
functions of bounded variation on $I$ (modulo the Lebesgue measure)
under the bounded variation norm. For $\phi$, $\psi \in BV$ and for
$n \geq 1$, we consider the correlation function
$$  C_{\phi, \psi}(n):= \int_I \phi(g^n(x)) \psi(x) d\mu(x) -
            \int_I \phi(x) d\mu(x)  \int_I \psi(x) d\mu(x). $$
\begin{definition} \label{def-mr}
Let $g$ be any piecewise expanding map on $I$. The {\em mixing rate}
$\tau(g)$ of $g$ is given by
$$ \tau(g):= \inf\{ \tau \geq 0\; : \;
   C_{\phi,\psi}(n):= O(\tau^n) \mbox{ as } n \rightarrow \infty \mbox{
     for all } \phi, \psi \in BV\}. $$
\end{definition}
If $g$ is topologically mixing then $\tau(g)<1$ by \cite[Corollary
  3.10]{Viana97}. Conversely, if $g$ is not topologically mixing then
$\tau(g)=1$, as can be seen be taking $\phi$, $\psi$ to be the
characteristic functions on suitable subsets of $I$.

We now restrict attention to maps of the form $g=\sigma \circ f$
with $f \in \F_m$ and $\sigma \in S_N$. In this situation, the
absolutely continuous invariant probability measure $\mu$ is always
the Lebesuge measure. We will consider the worst mixing rate
attained by the maps $\sigma \circ f$ as $\sigma$ varies.
\begin{definition}  \label{def-tmn}
Given $f \in \F_m$ and $N \geq m$, we set
$$   \tn(f):= \max\{ \tau(\sigma \circ f) : \sigma \in S_N \} . $$
\end{definition}

The main results of this paper will be expressed in terms of
$\tn(f)$. Our first result characterizes those $f$ such that $\sigma
\circ f$ is topologically mixing for all $\sigma \in S_N$.

\begin{theorem} \label{mixing-criterion}
Let $f \in \F_m$. Then $\tn(f)=1$ if and only if either
\begin{enumerate}
\item[(i)] $m \mid N$, or
\item[(ii)] $m \mid 2N$ and $f=\fzig$ or $\fizig$.
\end{enumerate}
In particular, if $m \nmid 2N$ then $\sigma \circ f$ is
topologically mixing for all $f \in \F_m$ and all $\sigma \in S_N$.
\end{theorem}

\begin{figure}[ht]
  \centering
  \includegraphics[width=4cm]{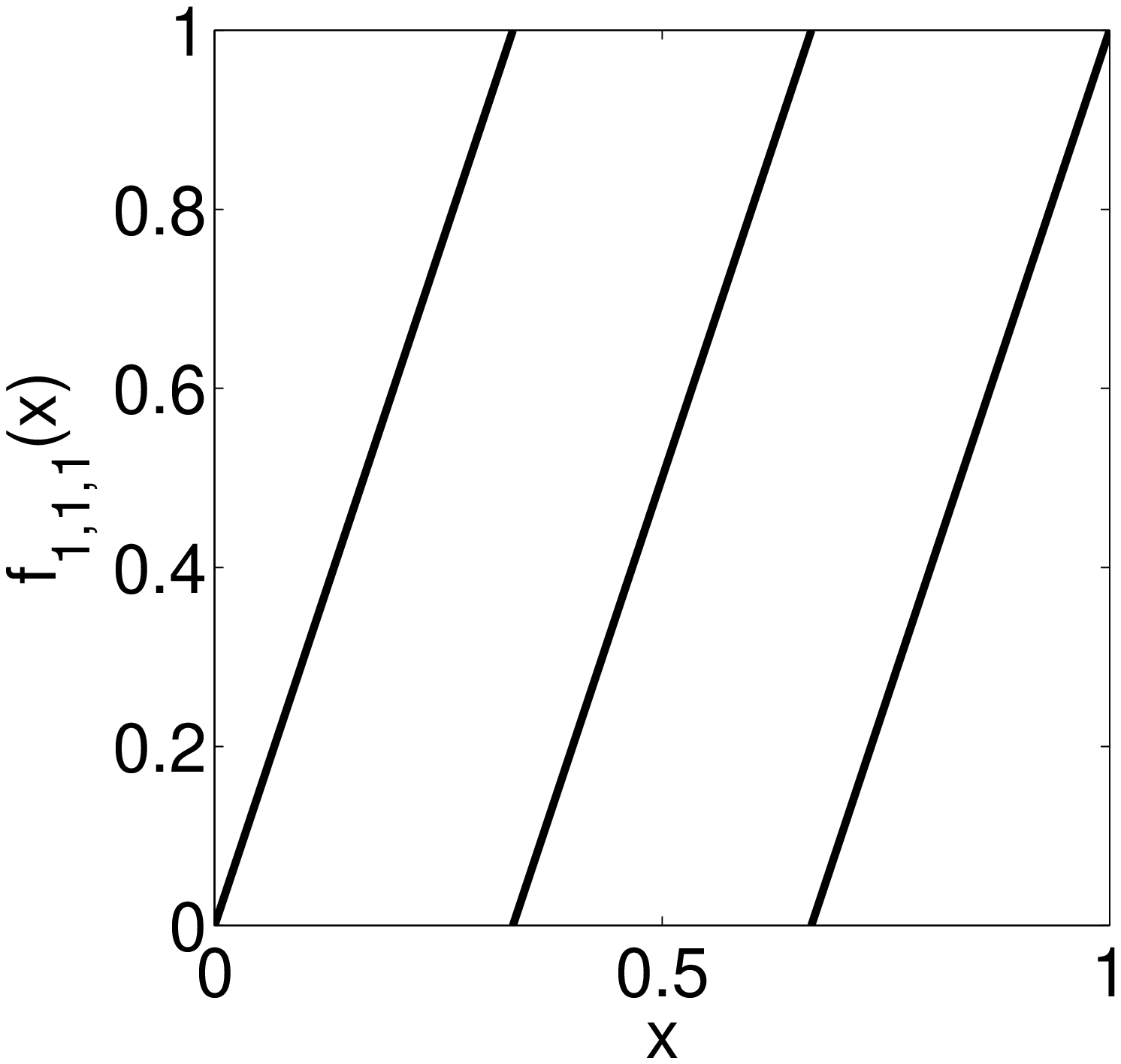}
  \includegraphics[width=4cm]{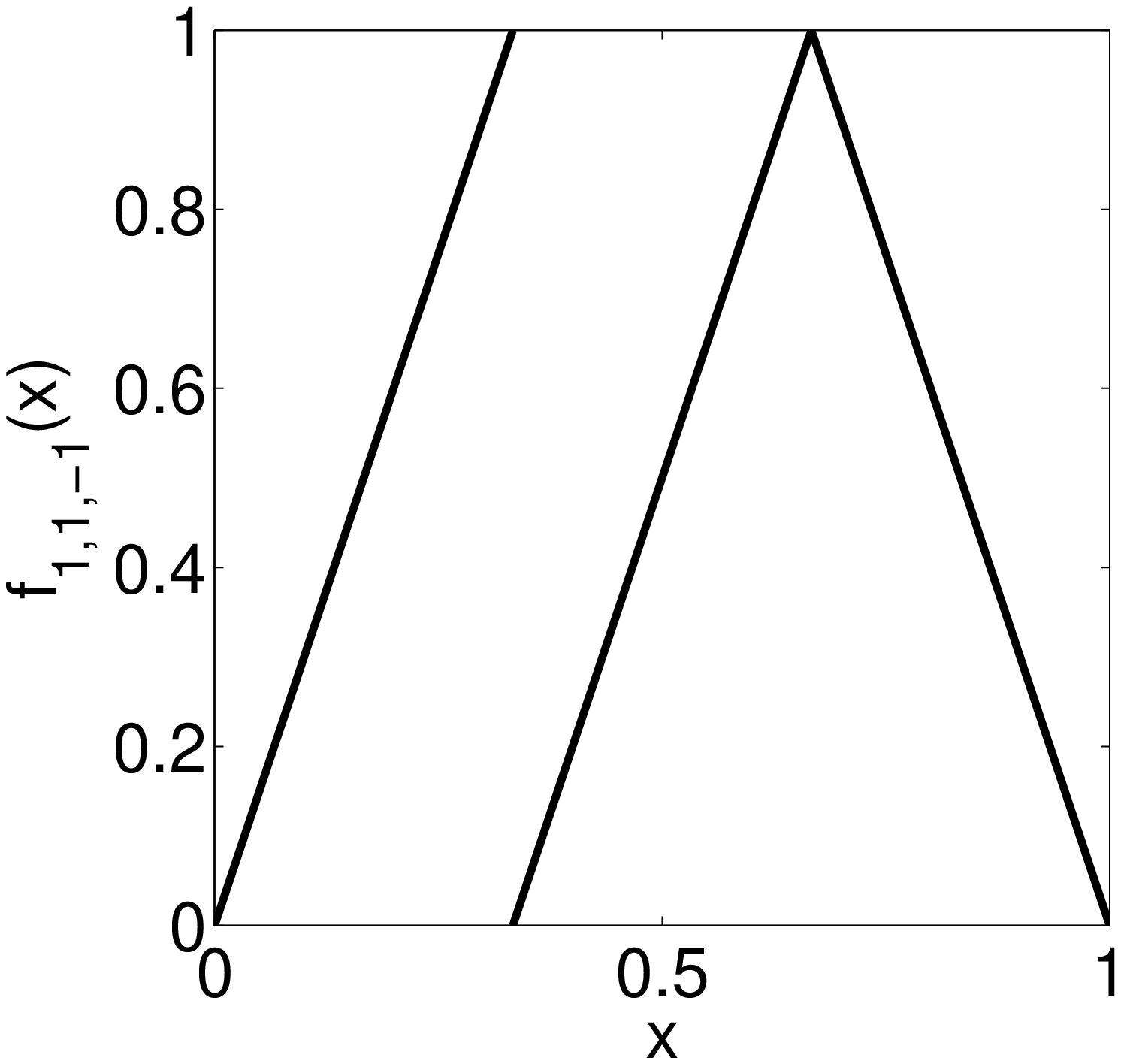}
  \includegraphics[width=4cm]{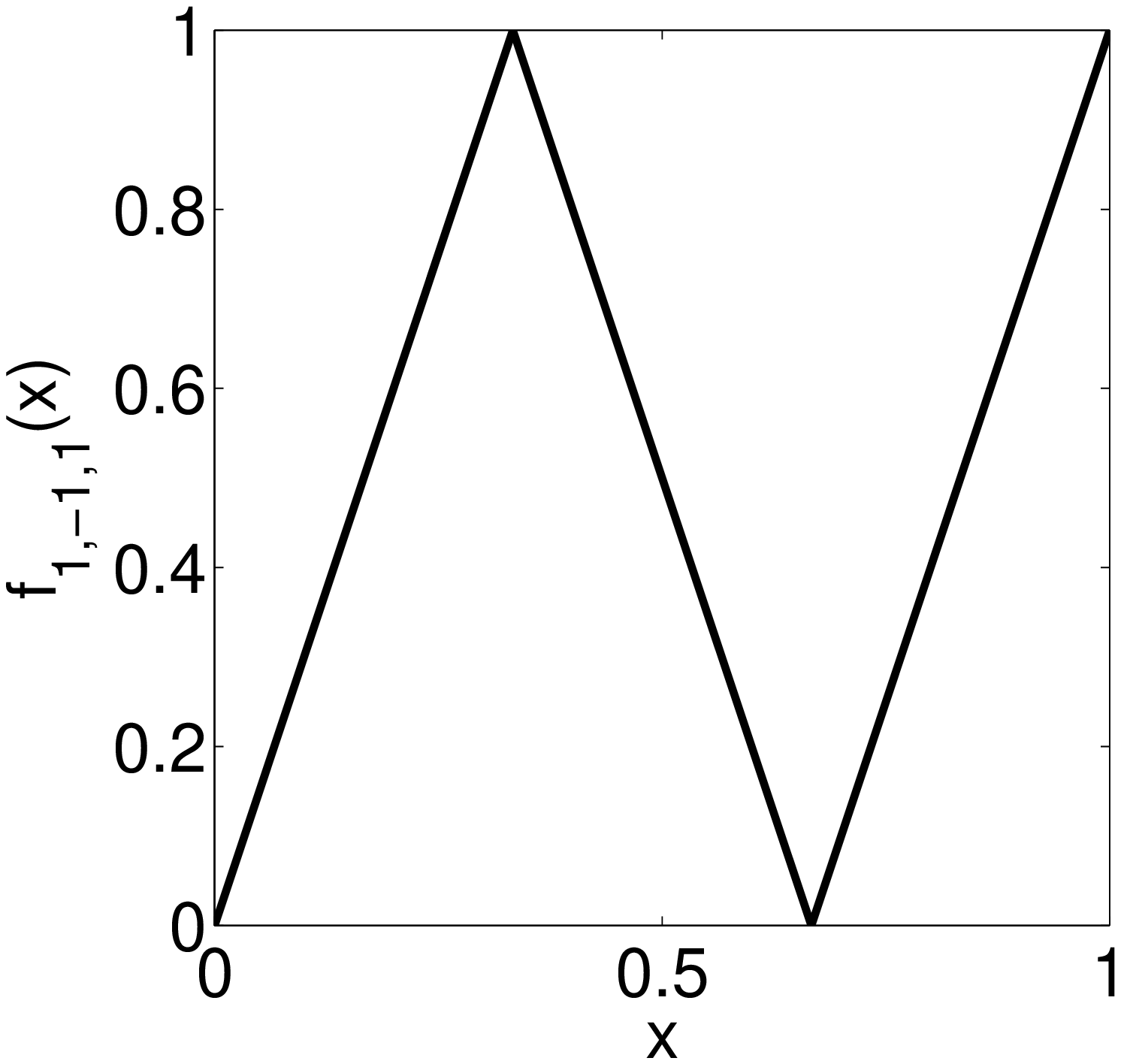}
  \caption{The family of maps $\mathcal{F}_{3}=\{f_{sf},f_{1,1,-1},f_{zz}\}$.}
  \label{fig_sfandzz}
\end{figure}

Although the family $\F_m$ contains $2^m$ maps, we will see in Lemma
\ref{symmetry} below that $\tn(f_{\epsilon_1,\ldots,
  \epsilon_m})$ is unchanged both by reversing the sequence
  $\epsilon_1, \ldots, \epsilon_m$ and by replacing $\epsilon_j$
  by $-\epsilon_j$ for all $j$. In particular
  $\tn(\fizig)=\tn(\fzig)$.
Up to these symmetries, $\F_2$ contains just two maps (namely, the
doubling map $\fsf$ and the tent map $\fzig=\ftent$), $\F_3$
contains three maps as illustrated in Figure \ref{fig_sfandzz},
$\F_4$ contains four maps, and $\F_5$ contains $10$ maps.

The maps $\fsf$ and $\fzig$ are extremal elements of $\F_m$ in that
they have respectively the smallest and largest number of changes of
slope. Our next result shows that, in a certain sense, no $f \in
\F_m$ can have worse asymptotic mixing behavior than $\fzig$.

\begin{theorem} \label{asymp-mixing}
For each odd $m \geq 3$ there is a constant $c(m)>0$ such that the
following holds. For any sequence of integers $N_i >m$ such that
$N_i \to \infty$ as $i \to \infty$ and $\gcd(m,N_{i})=1$ for each
$i$, and for any $f \in \F_m$, we have

\begin{equation} \label{mix-limit}
 \liminf_{i \to \infty} \frac{ 1 - \tni(f)}{1-\tni(\fzig)} \geq c(m).
\end{equation}
Explicitly, we may take $c(m) = 12/(m^4-m^2)$.

In the case $m=3$, we have the more precise statement:
\begin{equation}\label{equ_m=3}
\tn(f) \leq \tn(\fzig) \mbox{ for all } f \in \F_3 \mbox{ and for
  all } N \mbox{ with } \gcd(3,N)=1.
\end{equation}
\end{theorem}
Note that, by Theorem \ref{mixing-criterion}, the condition
$\gcd(m,N_i)=1$, with $m$ odd, ensures $\tni(f)$, $\tni(\fzig)<1$
for all $i$.

The proof of Theorem \ref{asymp-mixing} has two main ingredients.
One of these is a result of Fiedler \cite{Fiedler72} which bounds
the subleading eigenvalue of a doubly stochastic matrix away from
$1$. It is well-known that there is a strong connection between
mixing and spectral gaps for the transfer operator, but we are not
aware of any previous work where a general bound on the eigenvalues
of a stochastic matrix has been used to obtain explicitly
quantitative information on the limiting behavior of mixing rates.
The other ingredient is the following exact result on the zigzag
map.

\begin{theorem} \label{mr-zigzag}
Let $N \geq m  \geq 2$, and let $d=\gcd(m,2N)$. Then
\begin{equation} \label{zz-mr-gen}
 \tn(\fzig) = \frac{d \sin(m\pi/2N)}{m\sin(d\pi/2N)}.
\end{equation}
In particular, if $\gcd(m,2N)=1$ then
\begin{equation} \label{zz-mr-coprime}
    \tn(\fzig) = \frac{\sin(m\pi/2N)}{m\sin(\pi/2N)}.
 \end{equation}
\end{theorem}

The techniques used to prove Theorem \ref{mr-zigzag} are extensions
of those used for the stretch-and-fold map in \cite{BHZ}, and also
allow us to give the following improvement of \cite[Theorem
2.2(ii)]{BHZ} (where $m$ and $N$ were assumed coprime).

\begin{theorem} \label{sf-thm}
Let $N \geq m \geq 2$ and let $d=\gcd(m,N)$. Then
$$  \tn(\fsf) = \frac{d \sin(m\pi/N)}{m\sin(d\pi/N)}, $$
\end{theorem}

A case of particular interest is $\fzig$ for $m=2$, which is the
much-studied tent map. As $\gcd(2,2N)=2$, Theorem \ref{mr-zigzag}
merely tells us that $\tn(\ftent)=1$ for all $N$: there are some
permutations $\sigma \in S_N$ such that $\sigma \circ f$ is not
topologically mixing. Note, however, that iterates of $\ftent$ do
not share this property; the second iterate $\ftent^2$ of $\ftent$
is just the zigzag map for $m=4$, for which Theorem \ref{mr-zigzag}
gives $\tn(\fzig)<1$ for all odd $N \geq 5$.

In the cases where $\tn(f)=1$, it is natural to exclude those
permutations $\sigma$ where $\sigma \circ f$ is not topologically
mixing. We therefore give a modified version of Definition
\ref{def-tmn}:
\begin{definition} \label{def-thmn}
 $$  \tmix(f):= \max\{ \tau(\sigma \circ f) :
   \sigma \in S_N \mbox{ with } \sigma \circ f
                \mbox{ topologically mixing }\}  .  $$
\end{definition}
When $\sigma \circ f$ fails to be topologically mixing for some
$\sigma$, we have $\tmix(f) < \tn(f)=1$ and $\tmix(f)$ is obtained
by maximising over a proper subset of $S_N$. Although we have no
general method for approaching maximisation problems of this type,
we are able to give a partial result in the case of the tent map
$\ftent=\fzig$ for $m=2$. We give a lower bound for the worst mixing
rate $\tmix(\ftent)$. As explained in \S\ref{Markov}, the mixing
rate of $g= \sigma \circ f$ can be obtained as the subleading
eigenvalue of the transfer operator $\mathcal{L}_g$ of $g$. We are
able to obtain some geometric restrictions on the location in the
complex plane of the nonleading eigenvalues of $\mathcal{L}_g$ for
$g = \sigma \circ \ftent$.

\begin{figure}[ht]
  \centering
  \includegraphics[width=4cm]{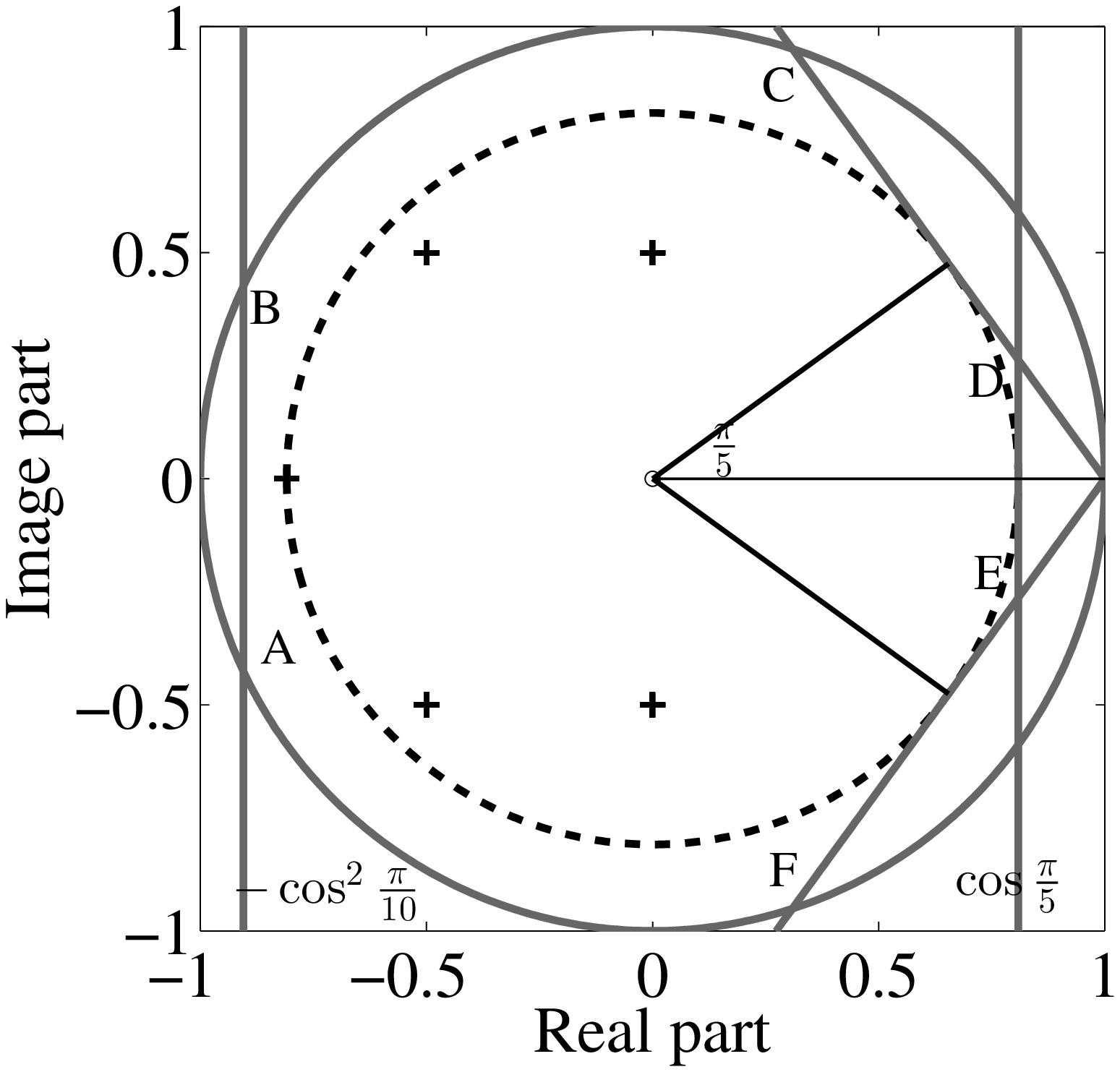}
  \includegraphics[width=4cm]{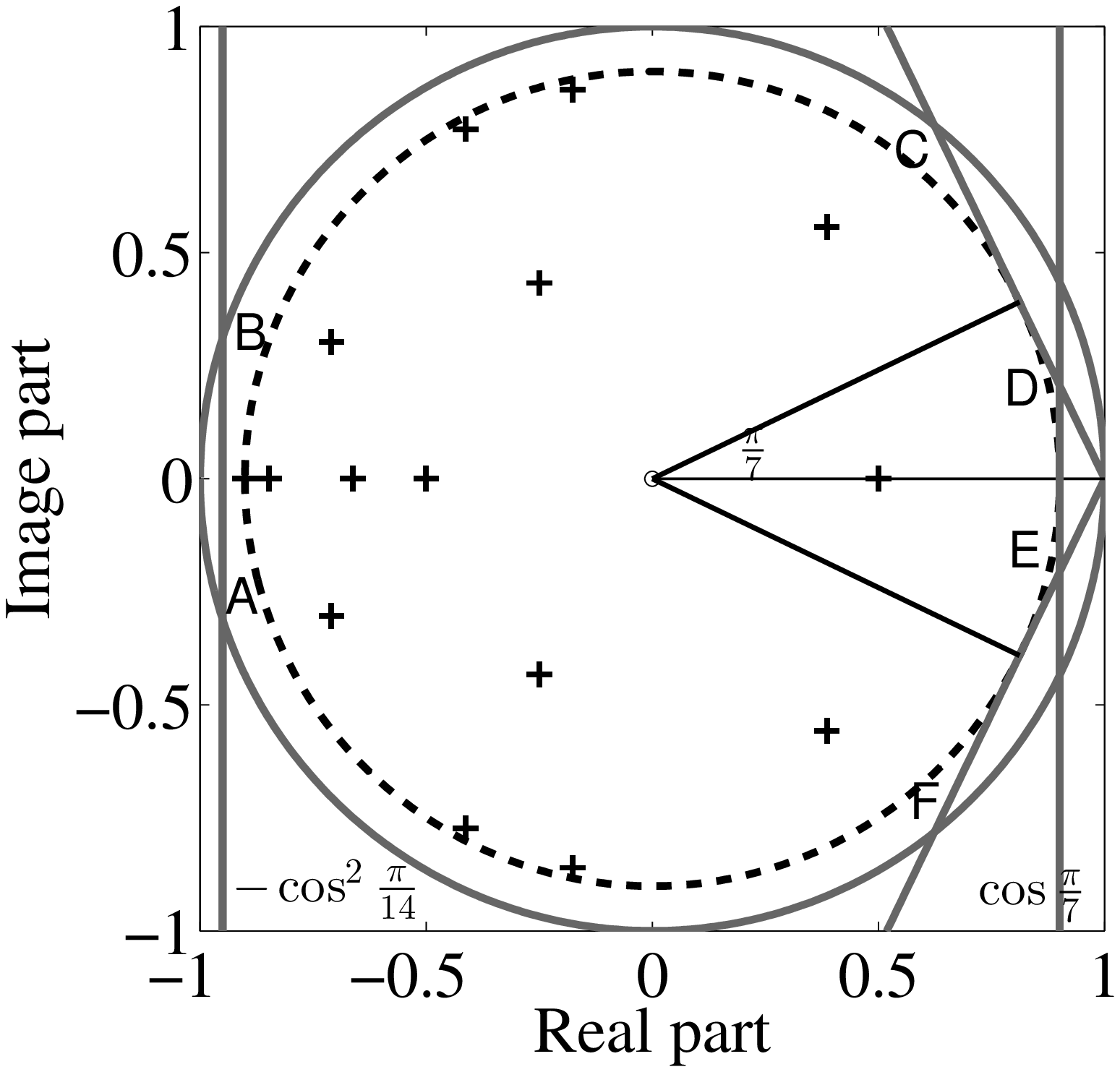}
  \includegraphics[width=4cm]{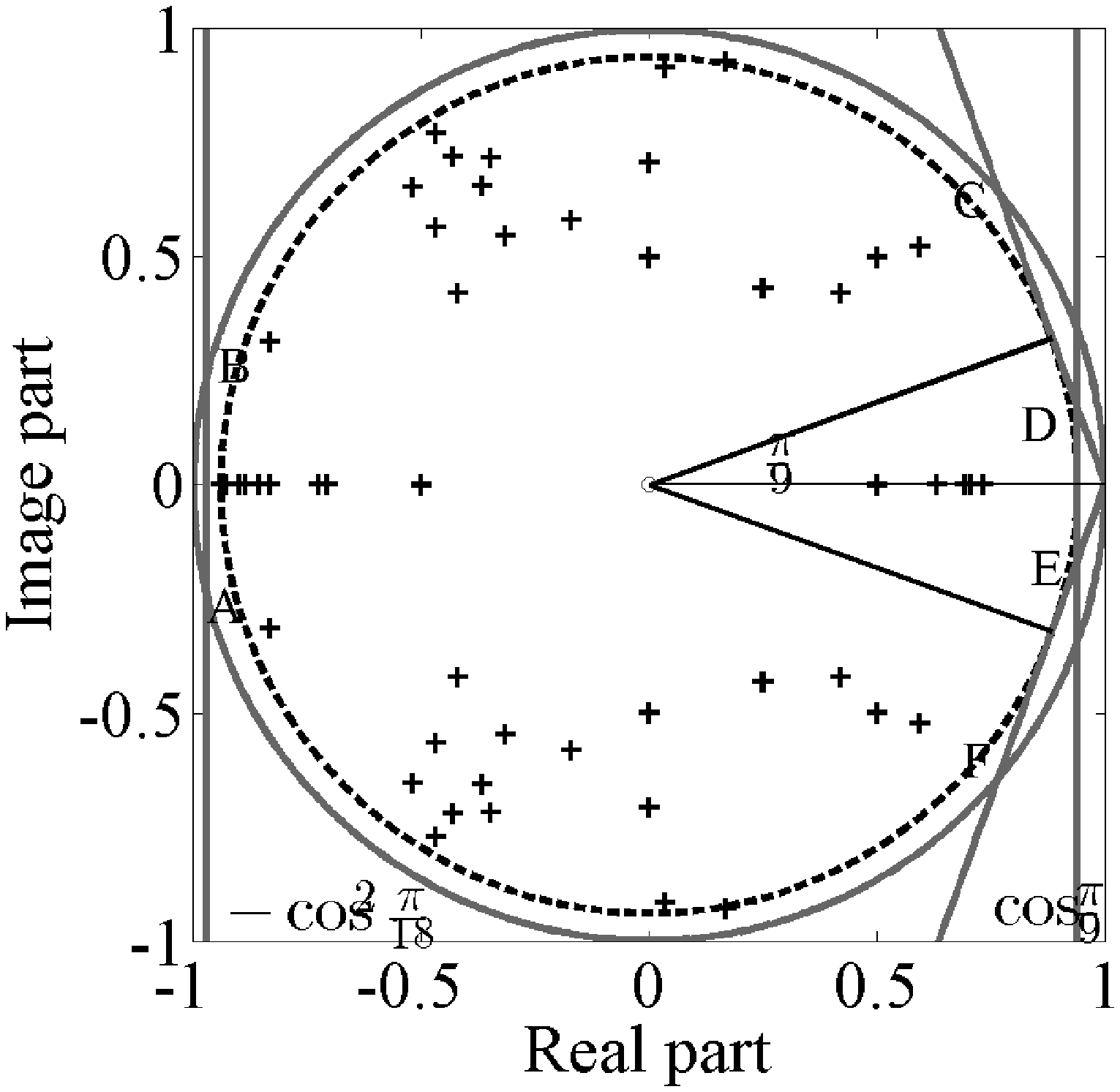}
  \caption{Location in the complex plane of the nonleading
    eigenvalues $\lambda$ of $\mathcal{L}_{\sigma\circ f_{tent}}$
    for $N=5,7,9$. Each such $\lambda$ is located in the convex subset
which is bounded \vspace{5pt} by the line segments $AB$, $CD$, $DE$,
$EF$ and the arcs $\overarc{BC}$, $\overarc{FA}$. The line segments
$CD$, $DE$, $EF$ are tangent to the disc centred at the origin with
radius $\cos(\frac{\pi}{N}).$}
  \label{Fig_eigenvalue}
\end{figure}

\begin{theorem} \label{tent}
Let $m=2$ and $N>2$ be odd. Then $\tmix(\ftent) \geq \cos(\pi/N)$.
Moreover, if $g=\sigma \circ \ftent$ for some $\sigma \in S_N$, and
$g$ is topologically mixing, then any nonleading eigenvalue
$\lambda$ of $\mathcal{L}_g$ lies in the compact convex region of
$\C$ given by the two inequalities
\begin{equation}\label{equ_spectral abscissa}
-\cos^{2}\left(\frac{\pi}{2N}\right)\leq\Re(\lambda)\leq\cos\left(\frac{\pi}{N}\right),
\end{equation}
and
\begin{equation}\label{equ_complex boundary}
\Re(\lambda)+|\Im(\lambda)|\tan\left(\frac{\pi}{N}\right)\leq1.
\end{equation}
\end{theorem}
Some examples of the regions given in Theorem \ref{tent} are
illustrated by Figure \ref{Fig_eigenvalue}.

Numerical calculations suggest that the lower bound on
$\tmix(\ftent)$ in Theorem~\ref{tent} is sharp, so that in fact
$\tmix(\ftent)= \cos(\pi/N)$.
\bigskip

\noindent {\bf Outline of the paper:}
\smallskip

In \S\ref{Markov}, we show how $\tn(f)$ can be calculated in terms
of Markov matrices for $f \in \F_m$. We then turn in \S\ref{exact}
to the exact values of $\tn(f)$ for the special cases $f=\fzig$ and
$\fsf$, giving the proofs of
Theorems~\ref{mr-zigzag}~and~\ref{sf-thm}. We also prove the
characterization of those $f$ with $\tn(f)=1$ in
Theorem~\ref{mixing-criterion}. Finally,
Theorems~\ref{asymp-mixing}~and~\ref{tent} are proved in
\S\ref{asymp}. We shall need a number of auxiliary results which
involve only linear algebra. To clarify the exposition, we postpone
their proofs to the Appendix~\S\ref{appendix}.

\section{Mixing rates and Markov matrices} \label{Markov}

In this section, we explain how mixing rates are related to the
transfer operator and can be calculated from certain matrices. We
restrict our discussion to maps of the form $g=\sigma \circ f$ with
$f \in \F_m$, $\sigma \in S_N$ and $N \geq m$, since these are the
only maps considered in this paper.

The transfer operator $\mathcal{L}_{g}:BV \to BV$ for $g$ is defined
by
\begin{equation}
\{\mathcal{L}_{g}\phi\}(x):=\sum_{g(y)=x}\frac{\phi(y)}{|g'|(y)},\quad\forall\phi\in
BV .
\end{equation}
The \emph{essential spectral radius} of $\mathcal{L}_g$ is
$$ r_{ess}(g):=\inf\{r\geq 0:
\lambda\in\operatorname{Spec}(\mathcal{L}_{g}), \,
|\lambda|>r\implies\lambda\,\textrm{is isolated}\}, $$ where the
spectrum is taken on the space $BV$. As $g$ has piecewise constant
slope $\pm m$, Keller's formula \cite[Theorem 1]{Keller84} gives
$r_{ess}(g)=1/m$.  The isolated eigenvalues $\lambda$ with
$|\lambda|>r_{ess}(g)$ are of finite multiplicity and satisfy
$|\lambda| \leq 1$. The eigenvalue $1$ always occurs. Moreover, $g$
is topologically mixing if and only if the eigenspace for the
eigenvalue $1$ has dimension $1$ and there are no other eigenvalues
$\lambda$ with
 $|\lambda|=1$. In this case, $\sup\left\{|\lambda|:\lambda\in
  \operatorname{Spec}(\mathcal{L}_{g})\backslash\{1\}\right\} =
  \tau(g)<1$, so the transfer operator has a spectral gap.

For $k \geq 1$, let $\mathcal{P}_k$ denote the partition of the unit
interval into the $k$ equal subintervals $I_j=[(j-1)/k, j/k]$ for $1
\leq j \leq k$. Each of our maps $g=\sigma \circ f$ is a Markov map
with respect to $\mathcal{P}_{Nm}$ with constant slope $\pm m$ on
each subinterval. Let $B(g,N)=(b_{ij})$ be the usual $\{0,1\}$
Markov transition matrix for $g$ on $\mathcal{P}_{Nm}$, so
$b_{ij}=1$ if and only if  $I_j^\circ \subset g(I_i)$. We note that
\begin{equation} \label{B-block}
     b_{i\;(h-1)m+k} = b_{i\;(j-1)m+1}  \mbox{ for } 1 \leq i \leq Nm,
     1 \leq j \leq N, 1 \leq k \leq m,
\end{equation}
since each subinterval of the partition $\mathcal{P}_{Nm}$ is mapped
by $g$ onto $m$ consecutive subintervals. The doubly stochastic
matrix $m^{-1} B(g,N)$ can be interpreted as the probability
transition matrix associated with $g$. This matrix is the Fredholm
matrix for $g$, as discussed in \cite[\S4.1]{BHZ}, from where we
have the following result.

\begin{lemma} \label{Phi-mr}
\begin{enumerate}
\item[(i)] The dynamical system given by $g$ is ergodic if and only if
$B(g,N)$ is irreducible, and is topologically mixing if and only if
$B(g,N)$ is primitive.
\item[(ii)] If $\lambda \in\C$ and $|\lambda|>r_{ess}(g)$ then
$\lambda\in\text{Spec}(\mathcal{L}_{g})$ if and only if $\lambda$ is
an eigenvalue of $m^{-1} B(g,N)$.
\end{enumerate}
\end{lemma}

We recall that a nonnegative matrix $M$ is said to be irreducible
(respectively, primitive) if, for all $i$, $j$, there is some $k
\geq 1$ with $m^{(k)}_{ij}>0$ (respectively, there is some $k \geq
1$ with $m^{(k)}_{ij}>0$ for all $i$, $j$), where
$M^k=(m^{(k)}_{ij})$ for $k \geq 1$.

Let $M$ be any nonnegative matrix with constant row and column sums
$c>0$.  By the Frobenius-Perron theorem, every eigenvalue $\lambda$
of $M$ satisfies $|\lambda| \leq c$. The {\em leading eigenvalue} of
$M$ is the eigenvalue $\lambda=c$, corresponding to the eigenvector
$(1, \ldots, 1)^T$. Moreover, the $(N-1)$-dimensional space
$\C^N_0=\{ (x_1, \ldots, x_N)^T \in \C^N : \sum_j x_j = 0 \}$ is
stable under $M$. We shall refer to eigenvalues of $M$ in $\C^N_0$
as {\em nonleading eigenvalues}. Note that for a nonleading
eigenvalue $\lambda$ we may have $\lambda=c$ (if $c$ is an
eigenvalue of $M$ of algebraic multiplicity $>1$), and we may have
$|\lambda|=c \neq \lambda$.  It is well-known that every nonleading
eigenvalue $\lambda$ of a stochastic matrix $M$ satisfies $\lambda
\neq 1$ (respectively, $|\lambda|<1$) if and only if $M$ is
irreducible (respectively, primitive).

It will be convenient to define an analogue of the mixing rate
$\tau(g)$ for matrices.
\begin{definition} \label{tau-matrix}
Let $M$ be a nonnegative $N \times N$ matrix with constant row and
column sums. If $N \geq2$, we define
$$  \tau(M) = \max \{ |\lambda| : \lambda \mbox{ is a nonleading
  eigenvalue of } M\}, $$
that is, $\tau(M)$ is the modulus of the subleading eigenvalue of
$M$.

In the degenerate case $N=1$, we define $\tau(M)=0$.
\end{definition}

We now explain how $\tau(g)$ can be obtained from an $N \times N$
matrix in place of the $Nm \times Nm$ matrix $B(g,N)$. Note that $f$
is, by definition, a Markov map on the partition $\mathcal{P}_m$,
but not in general on the partition $\mathcal{P}_N$. Thus $f$, and
hence also $g$, may change slope (from $m$ to $-m$ or vice versa) on
a subinterval of $\mathcal{P}_N$.
\begin{definition} \label{reduced-MM}
The {\em reduced Markov matrix} $A(g,N)$ is the $N \times N$ matrix
$(a_{ij})$ given by
$$ a_{ij} = \sum_{h=1}^m b_{(i-1)m+h\;(j-1)m+1} \in \{0, 1, 2\}. $$
\end{definition}

\begin{lemma} \label{small-mix-rate}
Let $g=\sigma \circ f$ with $f \in \F_m$ and $\sigma \in S_N$. Then
$$  \tau(B(g,N)) = \tau(A(g,N)).  $$
\end{lemma}
The proof of Lemma \ref{small-mix-rate} is given in the Appendix.

Since $\tau(g) \geq r_{ess}(g)=1/m$, the following result is then
immediate from Definition \ref{def-mr}, Lemma \ref{Phi-mr} and Lemma
\ref{small-mix-rate}:
\begin{corollary} \label{tau-phi}
$$  \tau(g)  =  \frac{1}{m} \max\left\{ 1,  \tau(A(g,N)) \right\}. $$
\end{corollary}

The mixing rates for the maps $f\in \F_m$ themselves (with no
permutation) are easily determined. In this case, $f$ is a Markov
map with respect to the partition $\mathcal{P}_m$ since each
subinterval of length $1/m$ is mapped homeomorphically to the whole
of $[0,1]$ (with either positive or negative orientation). We may
therefore apply Lemma \ref{Phi-mr} to the Fredholm matrix $F$ for
this partition (effectively taking $N=1$). But all the entries of
this $m \times m$ matrix are $1/m$, so its only nonleading
eigenvalue is $0$, occurring with multiplicity $m-1$. Thus
$\tau(F)=0$, and we have
\begin{equation} \label{mix-f}
  \tau(f)=1/m  \mbox{ for each } f \in \F_m.
\end{equation}

We next make explicit the relationship between the matrices $A(g,N)$
and $A(f,N)$ (respectively, $B(g,N)$ and $B(f,N)$) where $g=\sigma
\circ f$ with $\sigma \in S_N$. Let $P(\sigma)=(p_{ij})$ be the $N
\times N$ permutation matrix for $\sigma$:
\begin{equation} \label{perm-mx}
   p_{ij} = \begin{cases} 1 & \mbox{ if } j = \sigma(i), \cr
                                 0 & \mbox{otherwise,} \end{cases}
\end{equation}
and let $Q(\sigma)$ be the $Nm \times Nm$ permutation matrix
obtained by replacing each entry $1$ (respectively, $0$) in
$P(\sigma)$ by an $m \times m$ identity (respectively, zero) matrix.
Then
\begin{equation} \label{perm-AB}
  A(g,N)=A(f,N) P(\sigma), \qquad B(g,N)=B(f,N) Q(\sigma).
\end{equation}

It will be convenient to define an analogue for matrices of the
quantity $\tn(f)$ in Definition  \ref{def-tmn}.

\begin{definition} \label{tau-perm-matrix}
Let $M$ be a nonnegative $N \times N$ matrix with constant row and
column sums. Then
$$   \tperm(M) = \max \{ \tau(MP(\sigma)) : \sigma \in S_N \}. $$
\end{definition}

For any $\sigma \in S_N$, we have
\begin{equation} \label{tperm-change}
  \tperm(M P(\sigma) ) = \tperm(M) = \tperm(P(\sigma)M );
\end{equation}
the first inequality is immediate from the definition of $\tperm$,
and the second holds since $\C^N_0$ is stable under $P(\sigma)$ and
the conjugate matrices $MP(\sigma)$ and $P(\sigma)M$ have the same
eigenvalues.

It then follows from Definition \ref{def-tmn}, Corollary
\ref{tau-phi} and (\ref{perm-AB}) that
\begin{equation} \label{perm-mix-rate}
  \tn(f) =  \frac{1}{m}\; \max\left\{ 1, \tperm(A(f,N)) \right\}
  \mbox{ for all } f \in F_m.
\end{equation}

\begin{lemma}  \label{orthog-eigen}
Let $A$ be a nonnegative matrix with constant row and column sums.
Then $\tau(A A^T)=\tau(A^T A)\geq \tperm(A)^2$. If $A$ is symmetric
then $\tperm(A)=\tau(A)$.
\end{lemma}

The proof of Lemma \ref{orthog-eigen} is given in the Appendix. The
second assertion is essentially \cite[Lemma 4.3]{BHZ}.

We end this section by justifying the claim made in \S\ref{intro}
that $\tn(f)$ is unchanged by certain symmetries on $\F_m$.

\begin{lemma} \label{symmetry}
Let $(\epsilon_1, \ldots, \epsilon_m) \in \{\pm 1\}^m$, and let
$$ f=f_{\epsilon_1, \ldots, \epsilon_m}, \quad f_1=f_{-\epsilon_1,
  \ldots, -\epsilon_m}, f_2=f_{\epsilon_m, \ldots, \epsilon_1}. $$
Then $\tn(f_1)=\tn(f_2)=\tn(f)$.
\end{lemma}
{\bf Proof:} The matrix $A(f_1,N)$ is obtained by reversing the
order of the columns of $A(f,N)$, since each subinterval of length
$1/N$ is traversed in the opposite direction. Thus $A(f_1,N)=
A(f_1,N)P(\tau)$, where $\tau \in S_N$ is the permutation given by
$\tau(j)=N+1-j$ for $1 \leq j \leq N$. Similarly, $A(f_2,N)=
P(\tau)A(f,N)$. Hence by (\ref{tperm-change}),
$$ \tperm(A(f_1,N))=\tperm(A(f_2,N))=\tperm(A(f,N)). $$
and the result follows from (\ref{perm-mix-rate}). \hfill$\Box$

\section{Exact results}  \label{exact}

In this section, we prove Theorems \ref{mr-zigzag}, \ref{sf-thm} and
\ref{mixing-criterion}.

\subsection{Preliminary results}

\subsubsection{Eigenvectors of block matrices} \label{block-eigen}

Let $n \geq 1$ be a factor of $N$, say $N=dn$. We will relate $n
\times n$ matrices to certain $N \times N$ matrices, partitioned
into $d \times d$ blocks.

\begin{definition} \label{uparrow}
Let $A=(a_{ij})$ be an $n \times n$ matrix. Then $A^{\uparrow}$
denotes the $dn \times dn$ matrix obtained by replacing each entry
$a_{ij}$ of $A$ by a $d \times d$ matrix, each of whose entries is
$a_{ij}$.
\end{definition}

\begin{lemma} \label{mat-up}
Let $A$ be a nonnegative $d \times d$ matrix with constant row and
column sums. Then $\tau(A^{\uparrow}) = d \tau(A)$.
\end{lemma}

The proof of Lemma \ref{mat-up} is given in the Appendix.

\subsubsection{Matrices with a double symmetry property}

Let $J=J_N$ denote the ``backwards identity'' matrix
$(\delta_{i,N+1-j})_{1 \leq i,j \leq N}$, where $\delta_{ij}$ is the
Kronecker delta. Then $JM$ (respectively $MJ$) is the matrix
obtained from an $N \times N$ matrix $M$ by reversing the order of
its rows (respectively, columns). We also write $\tJ$ for the $2N
\times 2N$ matrix $J_{2N}$ of the same form.

We will require the following result, whose (easy) proof we leave to
the reader.

\begin{lemma} \label{sym-J-eigenvects-new}
Let $\tA$ be a $2N \times 2N$ matrix with constant column sums such
that $\tJ \tA = \tA = \tA \tJ$. Then
\begin{equation} \label{M-tM}
    \tA = \begin{pmatrix} A & AJ \\ JA & JAJ \end{pmatrix}
\end{equation}
for some $N \times N$ matrix $A$ with constant column sums.
Moreover,
$$ \tau(\tA) = 2 \tau(A). $$
\end{lemma}

\subsubsection{Circulant matrices} \label{circulant}

Let $m$, $N$ be natural numbers such that $1 \leq m \leq N$ and
$\gcd(m,N)=1$. We define
$$ \delta = \begin{cases} (1-m)/2 & \mbox{if $m$ is odd;} \cr
                                      (1-m+N)/2 & \mbox{if $m$ is even;} \end{cases}$$
and set $C=C(m,N)=(c_{ij})_{1 \leq i, j \leq N}$ with
$$  c_{ij} = \begin{cases} 1 & \mbox{if $j \equiv i + \delta + r \bmod N$
                                                    with $0 \leq r < m$;} \cr
                                       0 & \mbox{otherwise.} \end{cases} $$
Then $C$ is a symmetric circulant matrix. By \cite[(24) and
Proposition 6]{BHZ}, we have
\begin{equation} \label{max-circ-eval}
  \tau(C)  = \frac{\sin(m\pi/N)}{\sin(\pi/N)}.
\end{equation}

\begin{proposition} \label{D-sym}
Let $D=C+C J_N$. Then
$$ \tau(D) = \frac{2 \sin(m\pi/N)}{\sin(\pi/N)}. $$
\end{proposition}

The proof of Proposition \ref{D-sym} is given in the Appendix.

\subsection{Worst mixing rate for $\fsf$ and $\fzig$}

In this subsection, we prove Theorems~\ref{mr-zigzag}
and~\ref{sf-thm}, giving the worst mixing rate in the special cases
of the zigzag map $\fzig$ and the stretch-and-fold map $\fsf$. We
begin with the easier case $\fsf$.
\medskip

{\bf Proof of Theorem \ref{sf-thm}:} Let $d=\gcd(m,N)$, and write
$h=m/d$ and $n=N/d$. Let $A=A(\fsf,N)$. Then $A$ contains $n$
distinct rows, each occurring $d$ times. Let $C_0$ denote the
symmetric circulant matrix $C(h,n)$ as in \S\ref{circulant}, which
is defined since $\gcd(h,n)=1$, and let $C=C_0^{\uparrow}$ be the
corresponding $N \times N$ block matrix as given by Definition
\ref{uparrow}.  Then $A$ and $C$ have the same rows, so
$C=P(\sigma)A$ for some permutation $\sigma \in S_N$. Hence
$\tperm(A)=\tperm(C)$. Since $C$ is symmetric, $\tperm(C)=\tau(C)$
by Lemma \ref{orthog-eigen}.  But by Lemma \ref{mat-up} and
(\ref{max-circ-eval}), we have
$$ \tau(C) = d \tau(C_0) = \frac{d \sin(h\pi/n)}{\sin(\pi/n)}
   = \frac{d\sin(m\pi/N)}{\sin(d\pi/N)}  . $$
Then (\ref{perm-mix-rate}) gives
$$ \tn(\fsf) = \frac{1}{m} \tperm(A) =
    \frac{d\sin(m\pi/N)}{m \sin(d\pi/N)}  . $$
\hfill$\Box$

We now turn to the $m$-fold zigzag map $\fzig$.

{\bf Proof of Theorem \ref{mr-zigzag}:} We will relate the matrices
$A=A(\fzig,N)$ and $\tA=A(\fsf,2N)$. Explicitly, we have
$\tA=(\tilde{a}_{ij})_{1 \leq i,j  \leq 2N}$, where
$$ \tilde{a}_{ij} = \begin{cases} 1 & \mbox{if $j \equiv mi-r \bmod 2N$
    with $0 \leq r <m$;} \cr
     0 &  \mbox{otherwise.} \end{cases} $$
Let $J=J_N$ and $\tJ=J_{2N}$ as before. It is routine to verify that
$\tA$ commutes with $\tJ$. Write $E = \tA + \tA \tJ$. We have $\tJ E
= E = E \tJ$. The upper left-hand quarter of the $2N \times 2N$
matrix $E$ is precisely the $N \times N$ matrix $A(\fzig,N)$. (This
reflects the fact that we may regard $f$ as giving either a
dynamical system with states $1$, \ldots, $N$ corresponding the
subintervals of the partition $\mathcal{P}_N$, or a system with
states $1^+, \ldots N^+, N^-, \ldots, 1^-$; the states $j^+$, $j^-$
correspond to the subinterval $I_j$ with an orientation depending on
the slope of $f$. The construction of $E$ recombines states $j^+$,
$j^-$ into the single state $j$.) By Lemma
\ref{sym-J-eigenvects-new}, we then have
\begin{equation} \label{A2-doub}
     E = \begin{pmatrix} A & AJ \\ JA & J A J \end{pmatrix}  .
\end{equation}

Now consider the effect of replacing $\fzig$ by $g=\sigma \circ
\fzig$ for some $\sigma \in S_N$. This replaces $A$ by $AP(\sigma)$,
where $P(\sigma)$ is the permutation matrix corresponding to
$\sigma$, and so replaces $\tA$ by $\tA \tP(\sigma)$, where
$\tP(\sigma)$ is the $2N \times 2N$ permutation matrix
\begin{equation} \label{tP-sig}
   \tP(\sigma) =
  \begin{pmatrix} P(\sigma) & 0 \\ 0 & JP(\sigma)J \end{pmatrix}.
\end{equation}
Thus $E$ is replaced by the matrix
$$ \begin{pmatrix} A P(\sigma)  & A P(\sigma) J \\ J
  AP(\sigma)  & J A P(\sigma)  J \end{pmatrix}  = E \tP(\sigma). $$
Using Lemma \ref{sym-J-eigenvects-new} again, we have
$$ 2\tau(AP(\sigma)) = \tau(E\tP(\sigma)). $$

Next, let $d=\gcd(m,2N)$, and write $h=m/d$ and $n=2N/d$. Let
$C_0=C(h,n)$, and let $C=C_0^\uparrow$. Then $C=\hP(\pi)A(\fsf,2N)$
where $\hP(\pi)$ is the $2N \times 2N$ permutation matrix for some
$\pi \in S_{2N}$, as in (\ref{perm-mx}). Set $D=C+C\tJ$ and $D_0=C_0
+ C_0 J_n$, so that $D=D_0^\uparrow$, and let $\tP(\phi)$ be defined
as in (\ref{tP-sig}) with $\phi \in S_N$ given by
$$ \phi(i) = \begin{cases} \pi(i) & \mbox{ if } \pi(i) \leq N; \cr
                2N+1-\pi(i) & \mbox{ if } \pi(i) > N, \end{cases} $$
Then we have
$$ \hP(\pi) + \hP(\pi) \tJ = \tP(\phi) + \tP(\phi)
\tJ. $$ Since $\tA$ commutes with $\tJ$, it follows that
$D=\tP(\phi)E$. For each $\sigma \in S_N$, we therefore have
$$ 2 \tau(AP(\sigma)) = \tau(\tP(\phi)^{-1} D \tP(\sigma)) = \tau(D
\tP(\sigma \phi)). $$ Each matrix $\tP(\sigma \phi)$ is a
permutation matrix $\hP(\pi)$ for some $\pi \in S_{2N}$, but not
every $\pi \in S_{2N}$ can occur. Taking the maximum over $\sigma
\in S_N$, we have
$$ 2 \tperm(A) = \max_{\sigma\in S_N} \tau(D \tP(\sigma \phi))
               \leq \max_{\pi\in S_{2N}} \tau(D \hP(\pi))
               = \tperm(D). $$
On the other hand, taking $\sigma=\phi^{-1}$, we have
$$ 2 \tperm(A) \geq \tau(D). $$
But $D$ is a real symmetric matrix, so $\tperm(D)=\tau(D)$ by Lemma
\ref{orthog-eigen}. Hence we have $2 \tperm(A)=\tau(D)$.

We now calculate $\tau(D)$. Since $D=D_0^{\uparrow}$, we have
$\tau(D)= d \tau(D_0)$ by Lemma \ref{mat-up}. Since $D_0=C(h,n) +
C(h,n) J_n$, Proposition \ref{D-sym} tells us that
$$  \tau(D_0)  = \frac{2 \sin(h\pi/n)}{\sin(\pi/n)}
              = \frac{2 \sin(m\pi/2N)}{\sin(d\pi/2N)}, $$
so that
$$ \tau(D)  = \frac{2 d
  \sin(m\pi/2N)}{\sin(d\pi/2N)}.  $$

Finally, from (\ref{perm-mix-rate}), we have
$$ \tn(\fzig) = \frac{1}{m} \tperm(A) = \frac{1}{2m} \tau(D) =
 \frac{d \sin(m\pi/2N)}{m \sin(d\pi/2N)}. $$
This completes the proof of Theorem \ref{mr-zigzag}. \hfill$\Box$

\subsection{When the mixing rate is $1$}

In this subsection, we prove Theorem \ref{mixing-criterion}.  To do
so, we will use the facts that $\tn(\fsf)=1$ if $m \mid N$ and that
$\tn(\fzig)=1$ if $m \mid 2N$. These follow from
Theorem~\ref{sf-thm} (or \cite[Theorem 2.1(ii)]{BHZ}) and from
Theorem~\ref{mr-zigzag} respectively, but could easily be verified
directly.

Let $A=(a_{ij})$ be any nonnegative matrix with constant row and
column sums $m>0$. We introduce a corresponding ``row relation'' $R$
on the indexing set $\II=\{1,\ldots, N\}$ by writing $i \relR j$ if
and only if there is some $h \in \II$ with $a_{h i}a_{h j} \neq 0$.
Thus $i \relR j$ if and only if there is some state $h$ from which
both $i$ and $j$ can be reached in one step.  The relation $R$ is
reflexive and symmetric, but not necessarily transitive. We write
$\approx$ for the transitive closure of $R$. Then $\approx$ is an
equivalence relation on $\II$. Since both the row and column sums of
$A$ are constant, it follows that if $a_{h i} \neq 0$, then there is
a path from $i$ back to $h$. Thus if $i \relR j$ then $i$ and $j$
are in the same irreducible component of $\II$ for $A$. (The
converse need not to be true; for example, when $m=2$ the matrix
$A(\fzig,3)$ has $2$ equivalence classes under $\approx$ but only
one irreducible component.)

We now apply this to our reduced Markov matrices $A(f,N)$.

\begin{lemma} \label{irred-cond}
Let $f \in \F_m$ and $N \geq m$, and suppose that neither of the
conditions
\begin{enumerate}
\item[(i)] $m \mid N$, or
\item[(ii)] $m \mid 2N$ and $f=\fzig$ or $\fizig$.
\end{enumerate}
from Theorem \ref{mixing-criterion} hold. Then $A(f,N)$ is
irreducible, and $A(f,N)^T A(f,N)$ is primitive.
\end{lemma}
{\bf Proof:} Let $f=f_{\epsilon_1, \ldots, \epsilon_m}$ with
$\epsilon_1, \ldots, \epsilon_m = \pm1$, let $A=(a_{ij})$ be the
matrix $A(f,N)$, and let $\approx$ be the equivalence relation
corresponding to $A$ as above. To show that $A$ is irreducible, we
will verify that $\II$ consists of a single equivalence class under
$\approx$.

We suppose that $\epsilon_1=+1$, the case $\epsilon_1=-1$ being
similar.

As (i) does not hold, we may write $N=qm+r$ with $0 < r <m$. Now $f$
has slope $+m$ on $(0,1/m)$, so we have $a_{i,im-m+1}=a_{i,im-m+2}=
\cdots = a_{i,im}= 1$ for $1 \leq i \leq q$, and $a_{q+1,qm+1},
\ldots, a_{q+1,N}>0$. Thus the first $q$ rows of $A$ yield
\begin{equation} \label{approx1}
 im-m+1 \approx im-m+2 \approx \cdots \approx im \mbox{ for } 1 \leq i
 \leq q,
\end{equation}
and the next row yields
\begin{equation} \label{approx2}
  qm+1 \approx qm+2 \approx \ldots \approx N.
\end{equation}
This shows that there are at most $q+1$ equivalence classes under
$\approx$.

Now since neither (i) nor (ii) holds, some row of $A$ starts with a
block of $0$'s of length $<m$. Thus, for some $j$, we have $a_{j m}
a_{j\;m+1} \neq 0$. Then, for $2 \leq i \leq m$, we have $a_{k \;
im} a_{k-1 \; im+1}>0$ with $k=j+i-1$ or $j-i+1$. Hence $im \approx
im+1$ for $1 \leq i \leq m$. Together with (\ref{approx1}) and
(\ref{approx2}), this shows that there is a single equivalence class
under $\approx$. Hence $A$ is irreducible.

We now consider the matrix $B=A^T A$, which has constant row and
column sums $m^2$. Let us write $R_A$, $R_B$ for the row relations
$R$ corresponding to $A$, $B$ respectively, and similarly for
$\approx_A$, $\approx_B$. If $i \relR_A j$ then there is some $h$
with $a_{hi}$, $a_{hj}>0$. As the column sums of $A^T$ are $m>0$,
there is some $g \in \II$ such that $a^T_{gh}>0$. We then have
$b_{gi}$, $b_{gj}>0$, so $i \approx_B j$. As there is a single
equivalence class under $\approx_A$, the same must be true for
$\approx_B$. Thus $B$ is irreducible. This means that every
nonleading eigenvalue $\lambda$ of $B$ satisfies $\lambda \neq m^2$.
But $B$ is positive semidefinite symmetric matrix, so its
eigenvalues are real and nonnegative. Hence $|\lambda| \neq m^2$,
showing that $B$ is primitive. \hfill$\Box$

\medskip

{\bf Proof of Theorem \ref{mixing-criterion}:} If $\tn(f)=1$ then by
(\ref{perm-mix-rate}) there is some $\sigma \in S_N$ such that the
matrix $A(\sigma \circ f, N)=A(f,N)P(\sigma)$ is not irreducible.
Then $A(\sigma \circ f,N)$ has subleading eigenvalue $\lambda$ with
$|\lambda|=m$, so that $A(\sigma \circ f,N)^T A(\sigma \circ
f,N)=P(\sigma)^{-1} A(f,N)^T A(f,N) P(\sigma)$ has subleading
eigenvalue $m^2$, and hence so does its conjugate $A(f,N)^T A(f,N)$.
Thus $A(f,N)^T A(f,N)$ is not primitive, and by Lemma
\ref{irred-cond}, either (i) or (ii) must hold. This proves the
``only if'' direction of Theorem \ref{mixing-criterion}.

We now prove the converse. First, suppose that (i) holds, say
$N=qm$. Then $\tn(\fsf)=1$ by Theorem \ref{sf-thm}. Now for any $f
\in \F_m$, the the first $q$ rows of $A(f,N)$ are the same as the
corresponding rows in $A(\fsf,N)$ (possibly in reverse order), and
each subsequent block of $q$ rows simply repeats the first block
(possibly in reverse order). Thus $A(f,N)$ is obtained from
$A(\fsf,N)$ by some permutation of the rows, so that
$\tn(f)=\tn(\fsf)$. Hence $\tn(f)=1$ for all $f \in \F_m$. Finally,
suppose that (ii) holds. Then $\tn(\fzig)=1$ by Theorem
\ref{mr-zigzag}, and $\tn(\fizig)=\tn(\fzig)$ by Lemma
\ref{symmetry}.
\hfill$\Box$

\section{Asymptotic results}  \label{asymp}

In this section we prove Theorems \ref{asymp-mixing} and \ref{tent}.

\subsection{Eigenvalues of stochastic matrices}

There is a considerable literature on the eigenvalues of stochastic
matrices, see for instance \cite{Fiedler72, Fied75, FP75,
  John79, KS78, RT}. We review those results on this topic which
we will use.

We begin with Fiedler's bound \cite{Fiedler72}. Let $M=(m_{ij})$ be
a doubly stochastic $N \times N$ matrix. Fiedler defined the {\em
index of irreducibility} of $M$ to be
$$ \mu (M) = \min_{\mathcal{S}} \left\{ \sum_{i\in \mathcal{S}, j \not
  \in \mathcal{S}} m_{ij} \right\}, $$
where the minimum is over sets $\emptyset \subsetneq \mathcal{S}
\subsetneq \{1, \ldots, N\}$. Note that $0 \leq \mu(M) \leq 1$, with
$\mu(M)=0$ if and only if $M$ is reducible.

\begin{lemma} \cite[Theorem 3.4]{Fiedler72}  \label{Fiedler-bound}
Let $\lambda$ be any nonleading eigenvalue of the doubly stochastic
$N \times N$ matrix $M$. Then $|1 - \lambda| \geq 2 \bigl(1 - \cos
(\pi/N) \bigr) \mu(M)$.
\end{lemma}

A related result of Fiedler and Ptak bounds the eigenvalues away
from $-1$ in the case that $N$ is odd:

\begin{lemma} \cite[Theorem 3.4]{FP75} \label{F-Ptak}
Let $M$ be a doubly stochastic $N \times N$ matrix with $N$ odd.
Then any eigenvalue $\lambda$ of $M$ satisfies $|1 + \lambda| \geq
\bigl(1 - \cos (\pi/N) \bigr) \mu(M)$.
\end{lemma}

We next mention a result of Kellogg and Stephens which bounds the
eigenvalues of a nonnegative matrix in terms of the associated
directed graph. More precisely, given a nonnegative matrix
$G=(g_{ij})$ of order $N$, the directed graph $\mathcal{G}$ of $G$
has vertices $\{1,\cdots N\}$ and a directed edge $(i,j)$ if
$g_{ij}>0$. A sequence $i_{1},i_{2},\cdots,i_{k}$ of distinct
vertices is a \emph{closed circuit} of $\mathcal{G}$ if
$\mathcal{G}$ contains the edges
$(i_{1},i_{2}),(i_{2},i_{3}),\ldots, (i_{k},i_{1})$. Therefore, the
distinct vertices $i_{1},i_{2},\ldots,i_{k}$ form a circuit if and
only if the product $g_{i_{1}i_{2}}\cdots g_{i_{k},i_{1}}\neq0$. The
length of this circuit is $k$.

\begin{lemma} \cite[Theorem 1]{KS78} \label{KS}
Let $G$ be a nonnegative matrix with directed graph $\mathcal{G}$
and spectral radius $\rho$. Let $\kappa$ be the length of the
longest circuit of $\mathcal{G}$. If $\kappa=2$, all the eigenvalues
of $G$ are real. If $\kappa>2$, each eigenvalue $\lambda$ of $G$
satisfies
$$
\Re(\lambda)+|\Im(\lambda)|\tan\left(\frac{\pi}{\kappa}\right) \leq
\rho.
$$
\end{lemma}

We will also need the following standard result.

\begin{lemma} {\em (Weyl's inequalities)}  \cite[Theorem
      III.2.1]{Bhatia} \label{Weyl-ineq}
Let $A$, $B$ be two real symmetric $N \times N$ matrices with
eigenvalues $\alpha_1 \geq \ldots \geq \alpha_N$ and $\beta_1 \geq
\ldots \geq \beta_N$.  Let $\gamma_1 \geq \ldots \geq \gamma_N$ be
the eigenvectors of $A+B$. Then
$$ \gamma_j \leq \alpha_i + \beta_{j-i+1} \mbox{ for } i \leq j, \qquad
   \gamma_j \geq \alpha_i + \beta_{j-i+N} \mbox{ for } i \geq j. $$
\end{lemma}

\subsection{The extremal property of $\fzig$}  \label{extremal}
In this subsection, we prove Theorem \ref{asymp-mixing}, using
Theorem \ref{mr-zigzag} and Fiedler's bound, Lemma
\ref{Fiedler-bound}. We first prove \eqref{mix-limit}.

Fix $i$ and let $A=A(f,N_i)$. We consider the doubly stochastic
matrix $m^{-2} A^T A$. As $\gcd(m,N_i)=1$ and $m$ is odd, it follows
from Lemma \ref{irred-cond} that this matrix is primitive, and hence
irreducible. As its entries lie in $m^{-2} \Z$, we have
$\mu(m^{-2}A^T A) \geq m^{-2}$. Moreover, $m^{-2} A^T A$ is
symmetric and positive semidefinite, so its eigenvalues are real and
nonnegative. Thus, writing $\theta_i=\pi/2N_i$, Lemma
\ref{Fiedler-bound} gives
$$ 1- \tau(m^{-2} A^T A) \geq 2m^{-2} \bigl( 1 - \cos (2 \theta_i)
\bigr) = 4m^{-2} \sin^2 \theta_i. $$ But $\tau(m^{-2} A^T A) \geq
m^{-2} \tperm(A)^2$ by Lemma \ref{orthog-eigen}. Hence
\begin{eqnarray*}
 1 - m^{-1} \tperm(A) & = & \frac{ 1 - m^{-2} \tperm(A)^2}{ 1 + m^{-1}
   \tperm(A)} \\
                                  & \geq & \frac{1}{2} \left( 1 -
 \tau(m^{-2}A^T  A) \right) \\
                                   & \geq & 2 m^{-2} \sin^2  \theta_i.
\end{eqnarray*}

From (\ref{perm-mix-rate}), we have either $\tni(f)= m^{-1}
\tperm(A)$ or $\tni(f)= m^{-1}$. If $\tni(f) = m^{-1} \tperm(A)$, we
therefore have
$$ 1 - \tni(f) \geq 2m^{-2} \sin^2 \theta_i. $$
As $m \geq 3$, this also holds if $\tni(f) = m^{-1}$.

On the other hand, from Theorem \ref{mr-zigzag} we have
$$ 1-\tni(\fzig) = \frac{m \sin \theta_i - \sin(m\theta_i)}{m\sin \theta_i}. $$
Thus
$$ \liminf_{i \to \infty} \frac{1-\tni(f)}{1-\tni(\fzig)} \geq
\lim_{\theta \to 0}  \frac{2 m^{-1}\sin^3 \theta}{m\sin \theta
-\sin(m\theta)}
                           =   \frac{12}{m^4-m^2}. $$
This completes the proof of inequality \eqref{mix-limit}.

We now turn to the proof of (\ref{equ_m=3}). Thus we take $m=3$ and
suppose that $\gcd(3,N)=1$. By Lemma \ref{symmetry}, there are only
three maps $f \in \F_3$ to consider, namely $\fsf$, $\fzig$ and
$f_{1,1,-1}$. By Theorems \ref{mr-zigzag} and \ref{sf-thm}, we have
$\tn(\fsf) < \tn(\fzig)$, and it remains to show that $\tn(f) \leq
\tn(\fzig)$ for $f=f_{1,1,-1}$.

The $N \times N$ matrix $A(f,N)$ has constant row and column sums
$3$, its entries being $0$ or $1$ apart from one occurrence of $2$
in the final column. We define two related symmetric $\{0,1\}$
matrices $G=(g_{ij})_{1 \leq i,j  \leq N}$ and $P=(p_{ij})_{1 \leq
i,j \leq
  N}$ as follows:
$$ g_{ij} = \begin{cases} 1 & \mbox{if $i+j=N-1$ or $N$;} \cr
          1 & \mbox{if $(i,j)=(N-1,N)$ or $(N,N-1)$ or $(N,N)$;} \cr
            0 & \mbox{otherwise;} \end{cases} $$
and, writing $N=3k+r$ with $r \in \{1,2\}$,
$$ p_{ij} =  \begin{cases} 1 & \mbox{if $\lceil i/3 \rceil + \lceil
    j/3 \rceil = k+1$ and $i+j \equiv r+1 \bmod 3$}, \cr
     1 & \mbox{if $i=j>3k$;} \cr
            0 & \mbox{otherwise.} \end{cases} $$
Then the matrix $T=G+P$ is the unique symmetric matrix $T$ which can
be obtained by permuting the rows of $A(f,N)$. To estimate $\tn(f)$,
however, it is not sufficient to bound the eigenvalues of $T$ away
from $1$ using Lemma \ref{Fiedler-bound}, since $T$ is not in
general positive semidefinite. Instead, we use Weyl's inequalities
(Lemma \ref{Weyl-ineq}).

Let $\alpha_1 \geq \ldots \geq \alpha_N$ and $\beta_1 \geq \ldots
\geq \beta_N$ be the eigenvalues of the nonnegative real symmetric
matrices $G$, $P$ respectively. The rows of the anticirculant matrix
$G$ can be permuted to give a symmetric circulant matrix $C$, in
which each row has two consecutive entries $1$ and all other entries
$0$. We have $\tau(C)=2 \cos(\pi/N)$. (This follows from
(\ref{max-circ-eval}) when $N$ is odd, but is easily verified for
all $N$.) Thus, by Lemma \ref{orthog-eigen}, $\tperm(G) \leq 2 \cos
(\pi/N)$. Hence $\alpha_1= 2$ and $-2 \cos(\pi/N) \leq \alpha_N \leq
\alpha_2 \leq 2\cos(\pi/N)$. On the other hand,
 $-1 \leq \beta_N \leq \ldots \leq \beta_1 \leq 1$
since the only eigenvalues of the symmetric permutation matrix $P$
are $\pm 1$. Now let $\lambda_1=3 \geq \lambda_2 \geq \ldots \geq
\lambda_N$ be the eigenvalues of $T$. For $j \geq 2$, Lemma
\ref{Weyl-ineq} gives
$$ \lambda_j \leq \alpha_2 + \beta_{j-1} \leq 2 \cos (\pi/N)+1
  \mbox{ and }
  \lambda_j \geq \alpha_N + \beta_j \geq -2 \cos (\pi/N) -1. $$
Hence
$$ \tperm(T) \leq 1 + 2 \cos(\pi/N) = 3 - 4 \sin^2(\pi/2N) $$
so that
$$ \tn(f) \leq  \frac{1}{3} \bigl(3-4\sin^2(\pi/2N)\bigr). $$
On the other hand, by Theorem \ref{mr-zigzag},
$$ \tn(\fzig) = \frac{\sin(3\pi/2N)}{3 \sin(\pi/2N)} = \frac{1}{3}
\bigl(3-4\sin^2(\pi/2N)\bigr). $$ So we have $\tn(f) \leq
\tn(\fzig)$, as required.

\subsection{The tent map}
In this subsection, we prove our final theorem, which concerns the
tent map $\ftent$ for odd $N$.
\medskip

{\bf Proof of Theorem \ref{tent}:} Let $m=2$ and let $N=2s+1$. We
first show that $\tmix(\ftent) \geq \tn(\fsf)$. From Theorem
\ref{sf-thm}, we have
$$ \tn(\fsf) = \frac{\sin(2\pi/N)}{2\sin(\pi/N)} = \cos(\pi/N). $$
Let $A=A(\ftent,N)=(a_{ij})$; explicitly
$$ a_{ij} = \begin{cases} 1 & \mbox{if $i=h$ or $N+1-h$ and $j=2h-1$
    or $2h$ for $1 \leq h \leq s$;} \\
  2 & \mbox{if } i=s+1, j=N ; \\
  0 & \mbox{otherwise.} \end{cases} $$
We can permute the rows of $A$ to get the matrix $D=(d_{ij})$ such
that
$$ d_{ij} = \begin{cases} 1 & \mbox{if $i=1$ or $3$ and $j=1$ or $2$;} \\
   1 & \mbox{if $i=2h-2$ or $2h+1$ and $j=2h-1$ or $2h$
          for $2 \leq h \leq s$;} \\
  2 & \mbox{if } i=N-1, j=N ; \\
  0 & \mbox{otherwise.} \end{cases} $$
Then $D$ is irreducible since the associated directed graph contains
a circuit of length $N$ given by the sequence of edges
$$
   (1,2),(2,4),(4,6), \cdots, (N-3,N-1),(N-1,N),(N,N-2),\cdots,(5,3),(3,1).
$$
Thus $\tmix(\ftent) \geq \frac{1}{2}\tau(D)$.

The matrix $D$ has eigenvalue $0$ with multiplicity $s$ and
(leading) eigenvalue $2$ with multiplicity $1$. We claim that the
remaining $s$ nonleading eigenvalues are $\lambda_r=2 \cos(2\pi
r/N)$ for $1 \leq r \leq s$. We will then have $\tmix(\ftent) \geq
\frac{1}{m} \tau(D) = \frac{1}{2} \max_r |\lambda_r|=\cos(\pi/N)$,
as required.

To prove the claim, we observe that a {\em row} vector of the form
$$ (\alpha_s, \alpha_s, \alpha_{s-1}, \alpha_{s-1}, \ldots, \alpha_1,
\alpha_1, \alpha_0)$$ is  an eigenvector of $D$ with eigenvalue
$\lambda$ if and only if
$$ \alpha_s+\alpha_{s-1}=\lambda \alpha_s, \qquad
   \alpha_{h+1}+\alpha_{h-1} = \lambda \alpha_h \mbox{ for } 1 \leq h
   \leq s-1, \qquad 2\alpha_1=\lambda_0. $$
Fix $r$ with $1 \leq r \leq s$, let $\zeta=e^{2\pi i r/N}$, and set
$\alpha_j=\zeta^j+\zeta^{-j}$. Then the above equations hold with
$\lambda=\lambda_r$; this follows easily on noting that
$\alpha_0=2$, $\alpha_{h+1}+\alpha_{h-1}=\alpha_1 \alpha_h$ and
$\alpha_s=\alpha_{s+1}$. Thus $\lambda_r$ is indeed an eigenvalue of
$D$, as claimed.

Next we prove (\ref{equ_spectral abscissa}) and (\ref{equ_complex
  boundary}). Let $\sigma\in S_{N}$ be a
permutation such that $\sigma\circ f_{tent}$ is topologically
mixing, and let $U_{\sigma}=(u_{ij})=\frac{1}{2} A(\sigma\circ
f_{tent}, N)$ be the corresponding $N\times N$ primitive, doubly
stochastic matrix. Then $\mu(U_\sigma) \geq \frac{1}{2}$. Moreover,
for each $\emptyset\subsetneq \mathcal{S}\subsetneq \{1,\cdots,N\}$,
we have
\begin{equation}\label{equ_sum}
0<\sum_{i\in\mathcal{S},~j\notin\mathcal{S}}u_{ij} =
\sum_{i\notin\mathcal{S},~j\in\mathcal{S}}u_{ij}
=\sharp\mathcal{S}-\sum_{i\in\mathcal{S},~j\in\mathcal{S}}u_{ij}.
\end{equation}
Now let $V_{\sigma}=\frac{1}{2}(U_{\sigma}+U_{\sigma}^T)$. Then
$V_\sigma$ is symmetric and doubly stochastic, and it follows from
\eqref{equ_sum} that $\mu(V_{\sigma})=\mu(U_{\sigma})$. In
particular, $V_\sigma$ is irreducible.

Let the eigenvalues of the real symmetric stochastic matrix
$V_{\sigma}$ be
$1=\eta_{1}>\eta_{2}\geq\eta_{3}\geq\cdots\geq\eta_{N}>-1$. Let
$\yy$ be a unit eigenvector for a nonleading eigenvalue $\lambda$ of
$U_\sigma$. Then, writing $(\cdot, \cdot)$ for the usual complex
inner product on $\C^N$, we have
$(U_{\sigma}\yy,\yy)=\lambda(\yy,\yy)$ and
$(U_{\sigma}^T\yy,\yy)=(\yy,U_{\sigma}\yy)=\bar{\lambda}(\yy,\yy)$.
Hence,
$$
\Re(\lambda)=(V_{\sigma}\yy,\yy)\geq\eta_{N}.
$$
Furthermore, since $N$ is odd, it follows from Lemma \ref{F-Ptak}
that
$$ \Re(\lambda)\geq
 -1+\mu(V_{\sigma})\left( 1-\cos\left (\frac{\pi}{N} \right) \right)
 \geq-\cos\left(\frac{\pi}{2N}\right)^{2},
$$
giving the first inequality of \eqref{equ_spectral abscissa}.

On the other hand, if $\lambda$ is a nonleading eigenvalue of
$U_\sigma$ corresponding to the unit eigenvector $\yy$, then
$$ \eta_2=\max\limits_{||\xx||=1,\xx\in \C^N_0}(V_{\sigma}\xx,\xx)\geq
(V_{\sigma}\yy,\yy)=\Re(\lambda).
$$
Applying Fiedler's bound Lemma \ref{Fiedler-bound}, we have
$$
1-\eta_2\geq
2\left(1-\cos\left(\frac{\pi}{N}\right)\right)\mu(V_{\sigma})
 \geq 1-\cos\left(\frac{\pi}{N}\right),
$$
that is,
$$
\Re(\lambda)\leq \eta_2\leq\cos\left(\frac{\pi}{N}\right).
$$
This completes the proof of \eqref{equ_spectral abscissa}.

Finally, we prove the inequality \eqref{equ_complex boundary}. Let
$\kappa$ be the length of the longest circuit in the graph
associated to $U_\sigma$. As $U_\sigma$ is irreducible, we have
$\kappa>1$. Thus $2 \leq \kappa \leq N$, and, as $U_\sigma$ has
spectral radius $\rho=1$, the Kellogg-Stephens bound Lemma \ref{KS}
gives
$$ \Re(\lambda)+|\Im(\lambda)|\tan\left(\frac{\pi}{N}\right)\leq
 \Re(\lambda)+|\Im(\lambda)|\tan\left(\frac{\pi}{\kappa}\right)\leq
    1 $$
as required.

The bound $\kappa=N$ in the last part of the preceding proof is
attained by the matrix $D$ in the first part of the proof, since, as
noted above, this matrix contains a circuit of length $N$.
\hfill$\Box$

\section{Appendix: Proofs of linear algebra results}  \label{appendix}
We conclude by proving a few lemmas about linear algebra, which are
used in \S\ref{Markov} and \S\ref{exact}.
\subsection{Proof of Lemma \ref{orthog-eigen}}

Let $A$ be a real $N \times N$ matrix with eigenvalues $\lambda_1,
\ldots, \lambda_N$, numbered so that $|\lambda_1| \geq |\lambda_2|
\geq \ldots \geq |\lambda_N|$. The singular values of $A$ are the
real numbers $s_1 \geq s_2 \geq \ldots \geq s_N \geq 0$ such that
$s_1^2, \ldots, s_N^2$ are the eigenvalues of each of the symmetric,
positive semidefinite matrices $A A^T$ and $A^T A$ \cite[\S
I.2]{Bhatia}. For $1 \leq k \leq N$, followed by
\cite[(II.24)]{Bhatia}, we have
\begin{equation} \label{WH}
    \prod_{i=1}^k |\lambda_i| \leq \prod_{i=1}^k s_i;
\end{equation}

Now let $A$ also be nonnegative and have constant row and column
sums $m$. We may assume $m \neq 0$. Then $\lambda_1=s_1=m$,
$|\lambda_2|=\tau(A)$ and $s_2^2=\tau(A A^T) = \tau(A^T A)$. By
(\ref{WH}) with $k=2$, we have $m \tau(A) \leq m \sqrt{\tau(A^T
A)}$, so that $\tau(A)^2 \leq \tau(A A^T)$. Replacing $A$ by
$AP(\sigma)$ for $\sigma \in S_N$, we have $\tau(AP(\sigma))^2 \leq
\tau(AP(\sigma) P(\sigma)^T A^T)=\tau(AA^T)$. As this holds for all
$\sigma$, it follows that $\tperm(A)^2 \leq \tau(A A^T)$, giving the
first assertion of Lemma \ref{orthog-eigen}.

Suppose further that $A$ is symmetric. Then the eigenvalues
$\lambda_j$ are real and the corresponding eigenvalues $\ee_j$ can
be chosen to form an orthogonal basis of $\C^N$ with respect to the
usual complex inner product $( \cdot , \cdot)$. Then
$$ (A^T A \ee_i, \ee_j)= (A\ee_i, A\ee_j)=(\lambda_i \ee_i, \lambda_j \ee_j)
= \lambda_i \lambda_j (\ee_i, \ee_j)=\lambda_i^2 \delta_{ij}. $$
Thus the eigenvalues of $A^T A$ are the $\lambda_j^2$, so that
$\tau(A^T A) = \lambda_2^2 = \tau(A)^2$. As $\tperm(A) \geq \tau(A)
\geq 0$ by definition, and we have shown $\tau(A)^2=\tau(A^T A) \geq
\tperm(A)^2$, it follows that $\tperm(A)=\tau(A)$.

\subsection{Collapsing block matrices}

In this subsection, we prove Lemmas \ref{small-mix-rate} and
\ref{mat-up}. To do so, we need to ``collapse'' certain $N \times N$
matrices (viewed as $n \times n$ matrices of $d \times d$ blocks,
where $N=nd$), to obtain $n \times n$ matrices. For $1 \leq p \leq
N$, we write $p=(i-1)d+r$ with $1 \leq i \leq n$ and $1\leq r \leq
d$. Thus a subscript $p$ relates to the $r$th position in the $i$th
block. To emphasize this indexing by pairs $(i,r)$, we will write
the entries $b_{pq}$ of an $N \times N$ matrix $B$ as
$b_{ij}^{(rs)}$ in place of $b_{(i-1)d+r\:(j-1)d+s}$.

 We will say that $B$ has the {\em column
  block property for} $d$ if the $d$ columns in each block are
identical:
$$ b_{ij}^{(rs)} = b_{ij}^{(r1)} \mbox{ for } 1 \leq i, j \leq n
\mbox{ and } 1 \leq r, s \leq d. $$ We then define $B^\downarrow =(
b^\downarrow_{ij})$ to be the $n \times n$ matrix obtained by
replacing each block in $B$ by the sum of any of its columns:
$$ b^\downarrow_{ij} = \sum_{r=1}^{d} b_{ij}^{(rs)}. $$
(This sum is independent of $s$ because of the column block
property).

The next result is a slightly modified version of \cite[Lemma
4.2]{BHZ}.

\begin{lemma} \label{collapse}
Let $B$ be an $N \times N$ matrix with constant row and column sums,
and with the column block property for $d$. Then
$$ \tau(B^\downarrow) = \tau(B). $$
\end{lemma}
{\bf Proof:} For $1 \leq j \leq n$ and $1 \leq s \leq d$, let
$\vv_j^{(s)} \in \C^N$ be the vector whose component in position
$(i-1)d+r$ is as follows: if $s=1$ then
$$ v_{ji}^{(sr)}  = \begin{cases} 1 & \mbox{if } r=1, i=j, \cr
                                  0 & \mbox{otherwise}, \end{cases} $$
and if $s \neq 1$ then
$$ v_{ji}^{(sr)}  = \begin{cases} -1 & \mbox{if } r=1, i=j, \cr
                                  1 & \mbox{if } r=s, i=j, \cr
                                  0 & \mbox{otherwise}. \end{cases} $$
As $B$ has the column block property, we have $B\vv_j^{(s)}=\zero$
if $s \neq 1$. Thus the linear endomorphism $\theta$ of $\C^N$ given
by $B$ vanishes on the $n(d-1)$-dimensional subspace $W$ spanned by
the $\vv_j^{(s)}$ with $s \neq 1$. It follows that $\theta$ induces
an endomorphism $\overline{\theta}$ of the quotient space $\C^N/W$.
The $n$ cosets $\vv_j^{(1)}+W$ form a basis for this quotient space,
and the matrix for $\overline{\theta}$ with respect to this basis is
given by replacing each $d \times d$ block in $B$ by the sum of one
of its (identical) columns. This matrix is precisely $B^\downarrow$.
The eigenvalues of $\theta$ are those of $\overline{\theta}$, with
the same multiplicities, together with the eigenvalue $0$ with a
(further) multiplicity $n(d-1)$. Thus, if $n \geq 2$, we have
$\tau(B^\downarrow)= \tau(B)$. In the degenerate case $n=1$, we have
$\tau(B)=0$ since the only nonleading eigenvalue of $B$ is $0$ (with
multiplicity $d-1$), and by definition $\tau(B^\downarrow)=0$.
\hfill$\Box$

\medskip

{\bf Proof Lemma \ref{small-mix-rate}:} Let $g = \sigma \circ f$
with $\sigma \in S_N$. Then the $Nm \times Nm$ matrix $B(f,N)$ has
the column block property with $d=m$ by (\ref{B-block}), and hence
so does $B(g,N)$ since this property is preserved under multiplying
on the right by the block permutation matrix $Q(\sigma)$. Indeed, if
we partition $B(g,N)$ into $m \times m$ blocks, and collapse each
block to the sum of any one of its columns, the resulting $N \times
N$ matrix $B(g,N)^\downarrow$ is precisely $A(g,N)$, since each
block of $m$ consecutive rows corresponds to a single interval of
length $1/N$. Thus $\tau(A(g,N))=\tau(B(g,N))$ by Lemma
\ref{collapse}. \hfill$\Box$

\medskip

{\bf Proof of Lemma \ref{mat-up}:} Let $A$ be a nonnegative $d
\times d$ matrix with constant row and columns sums, and let
$B=A^{\uparrow}$ be the corresponding $N \times N$ matrix as in
Definition \ref{uparrow}. Then $B$ has the column block property,
and $B^\downarrow= d A$. Thus, using Lemma \ref{collapse}, $
\tau(A^{\uparrow}) = \tau(B) = \tau(B^\downarrow) = d \tau(A). $
\hfill$\Box$

\subsection{Proof of Proposition \ref{D-sym}}

Let $C=C(m,N)$ with $\gcd(m,N)=1$ be the symmetric circulant matrix
defined in \S\ref{circulant}. Let $\eta=e^{\pi i /N}$, let
$\omega_j=\eta^{2j}$ for all $j \in \Z$, and set
$$  \vv_j = (1, \omega_j, \omega_j^2, \ldots ,\omega_j^{N-1})^T \in
\C^N. $$ A routine calculation (cf.~\cite[p.~3394]{BHZ}) then shows
that $C \vv_j = \lambda_j \vv_j$ where
$$ \lambda_j = \begin{cases}
      m & \mbox{if } N \mid j, \\
     \displaystyle{ (-1)^{(m-1)j} \frac{\sin(mj\pi/N)}{\sin(j\pi/N)}}
       & \mbox{otherwise.}
              \end{cases}  $$
Now set $\ww_j = \eta^j \vv_j$. Then $J \ww_j =
\ww_{-j}=\ww_{2N-j}=-\ww_{N-j}$. For $1 \leq j < N/2$, define
$$  \ff_j  = \ww_j+\ww_{-j} = \ww_j - \ww_{N-j}~~\mbox{and}~~\gg_j  =  \ww_j-\ww_{-j} = \ww_j + \ww_{N-j}. $$ These vectors,
together with $\ww_0$ and (in the case that $N$ is even) $\ww_{N/2}$
clearly form a basis of $\C^N$. We will show that they are all
eigenvectors for $D=D(m,N)$. Certainly $D\ww_0=2m \ww_0$, giving the
leading eigenvalue $2m$. Now we calculate
\begin{eqnarray*}
 D \ff_j & = & \left( C \ww_j + CJ \ww_j - C \ww_{n-j} - CJ \ww_{n-j} \right) \\
 & = & \left( C \ww_j - C \ww_{n-j} - C \ww_{n-j} + C \ww_j \right) \\
 & = & 2 \left( \lambda_j \ww_j - \lambda_{n-j} \ww_{n-j}  \right) \\
 & = & 2 \lambda_j \ff_j.
\end{eqnarray*}
Thus $\ff_j$ is an eigenvector of $D$ with eigenvalue $2\lambda_j$.
Similar calculations show that $D \gg_j=\zero$ and (if $n$ is even)
$D \ww_{n/2}=\zero$. Hence the nonzero nonleading eigenvalues of $D$
are $2 \lambda_j$ for $1 \leq j < n/2$. The result now follows from
(\ref{max-circ-eval}).

\section*{Acknowledgments}
The authors would thank the "Probability in Dynamics" workshop at
IM-UFRJ, Rio de Janeiro, in 26-30 May, 2014. We also thank an
anonymous referee for his/her very helpful comments. Y. Zhang is
partially supported by FONDECYT, grant No. 3130622.

\end{document}